# PIECEWISE LINEAR REGULARIZED SOLUTION PATHS

By Saharon Rosset and Ji Zhu[1]

*IBM T. J. Watson Research Center and University of Michigan*

We consider the generic regularized optimization problem $\hat{\beta}(\lambda) = \arg\min_\beta L(\mathbf{y}, X\beta) + \lambda J(\beta)$. Efron, Hastie, Johnstone and Tibshirani [*Ann. Statist.* **32** (2004) 407–499] have shown that for the LASSO—that is, if $L$ is squared error loss and $J(\beta) = \|\beta\|_1$ is the $\ell_1$ norm of $\beta$—the optimal coefficient path is piecewise linear, that is, $\partial\hat{\beta}(\lambda)/\partial\lambda$ is piecewise constant. We derive a general characterization of the properties of (loss $L$, penalty $J$) pairs which give piecewise linear coefficient paths. Such pairs allow for efficient generation of the full regularized coefficient paths. We investigate the nature of efficient path following algorithms which arise. We use our results to suggest robust versions of the LASSO for regression and classification, and to develop new, efficient algorithms for existing problems in the literature, including Mammen and van de Geer's locally adaptive regression splines.

**1. Introduction.** Regularization is an essential component in modern data analysis, in particular when the number of predictors is large, possibly larger than the number of observations, and nonregularized fitting is likely to give badly over-fitted and useless models.

In this paper we consider the generic regularized optimization problem. The inputs we have are:

- A training data sample $X = (\mathbf{x}_1, \ldots, \mathbf{x}_n)^\top, \mathbf{y} = (y_1, \ldots, y_n)^\top$, where $\mathbf{x}_i \in \mathbb{R}^p$ and $y_i \in \mathbb{R}$ for regression, $y_i \in \{\pm 1\}$ for two-class classification.
- A convex nonnegative loss functional $L : \mathbb{R}^n \times \mathbb{R}^n \to \mathbb{R}$.
- A convex nonnegative penalty functional $J : \mathbb{R}^p \to \mathbb{R}$, with $J(0) = 0$. We will almost exclusively use $J(\beta) = \|\beta\|_q$ in this paper, that is, penalization of the $\ell_q$ norm of the coefficient vector.

Received September 2003; revised June 2006.

[1]Supported by NSF Grant DMS-05-05432.

*AMS 2000 subject classifications.* Primary 62J07; secondary 62F35, 62G08, 62H30.

*Key words and phrases.* $\ell_1$-norm penalty, polynomial splines, regularization, solution paths, sparsity, total variation.

We want to find

$$(1) \qquad \hat{\beta}(\lambda) = \underset{\beta \in \mathbb{R}^p}{\arg\min} \, L(\mathbf{y}, X\beta) + \lambda J(\beta),$$

where $\lambda \geq 0$ is the regularization parameter; $\lambda = 0$ corresponds to no regularization, while $\lim_{\lambda \to \infty} \hat{\beta}(\lambda) = 0$.

Many of the commonly used methods for data mining, machine learning and statistical modeling can be described as exact or approximate regularized optimization approaches. The obvious examples from the statistics literature are explicit regularized linear regression approaches, such as ridge regression [9] and the LASSO [17]. Both of these use squared error loss, but they differ in the penalty they impose on the coefficient vector $\beta$:

$$(2) \qquad \text{Ridge}: \qquad \hat{\beta}(\lambda) = \min_{\beta} \sum_{i=1}^{n} (y_i - \mathbf{x}_i^\top \beta)^2 + \lambda \|\beta\|_2^2,$$

$$(3) \qquad \text{LASSO}: \qquad \hat{\beta}(\lambda) = \min_{\beta} \sum_{i=1}^{n} (y_i - \mathbf{x}_i^\top \beta)^2 + \lambda \|\beta\|_1.$$

Another example from the statistics literature is the penalized logistic regression model for classification, which is widely used in medical decision and credit scoring models:

$$\hat{\beta}(\lambda) = \min_{\beta} \sum_{i=1}^{n} \log(1 + e^{-y_i \mathbf{x}_i^\top \beta}) + \lambda \|\beta\|_2^2.$$

Many "modern" methods for machine learning and signal processing can also be cast in the framework of regularized optimization. For example, the regularized support vector machine [20] uses the hinge loss function and the $\ell_2$-norm penalty:

$$(4) \qquad \hat{\beta}(\lambda) = \min_{\beta} \sum_{i=1}^{n} (1 - y_i \mathbf{x}_i^\top \beta)_+ + \lambda \|\beta\|_2^2,$$

where $(\cdot)_+$ is the positive part of the argument. Boosting [6] is a popular and highly successful method for iteratively building an additive model from a dictionary of "weak learners." In [15] we have shown that the AdaBoost algorithm *approximately* follows the path of the $\ell_1$-regularized solutions to the exponential loss function $e^{-yf}$ as the regularizing parameter $\lambda$ decreases.

In this paper, we concentrate our attention on (loss $L$, penalty $J$) pairings where the optimal path $\hat{\beta}(\lambda)$ is *piecewise linear* as a function of $\lambda$, that is, $\exists \lambda_0 = 0 < \lambda_1 < \cdots < \lambda_m = \infty$ and $\gamma_0, \gamma_1, \ldots, \gamma_{m-1} \in \mathbb{R}^p$ such that $\hat{\beta}(\lambda) = \hat{\beta}(\lambda_k) + (\lambda - \lambda_k)\gamma_k$ for $\lambda_k \leq \lambda \leq \lambda_{k+1}$. Such models are attractive because they allow us to generate the whole regularized path $\hat{\beta}(\lambda), 0 \leq \lambda \leq \infty$, simply by sequentially calculating the "step sizes" between each two consecutive $\lambda$ values and the "directions" $\gamma_1, \ldots, \gamma_{m-1}$. Our discussion will concentrate on $(L, J)$ pairs which allow efficient generation of the whole path and give statistically useful modeling tools.

A canonical example is the LASSO (3). Recently [3] has shown that the piecewise linear coefficient paths property holds for the LASSO, and suggested the LAR–LASSO algorithm which takes advantage of it. Similar algorithms were suggested for the LASSO in [14] and for total-variation penalized squared error loss in [13]. We have extended some path-following ideas to versions of the regularized support vector machine [7, 21].

In this paper, we systematically investigate the usefulness of piecewise linear solution paths. We aim to combine efficient computational methods based on piecewise linear paths and statistical considerations in suggesting new algorithms for existing regularized problems and in defining new regularized problems. We tackle three main questions:

1. What are the "families" of regularized problems that have the piecewise linear property? The general answer to this question is that the loss $L$ has to be a *piecewise quadratic* function and the penalty $J$ has to be a *piecewise linear* function. We give some details and survey the resulting "piecewise linear toolbox" in Section 2.

2. For what members of these families can we design efficient algorithms, either in the spirit of the LAR–LASSO algorithm or using different approaches? Our main focus in this paper is on direct extensions of LAR–LASSO to "almost-quadratic" loss functions (Section 3) and to nonparametric regression (Section 4). We briefly discuss some non-LAR type results for $\ell_1$ loss in Section 5.

3. Out of the regularized problems we can thus solve efficiently, which ones are of statistical interest? This can be for two distinct reasons:

(a) Regularized problems that are widely studied and used are obviously of interest, if we can offer new, efficient algorithms for solving them. In this paper we discuss in this context locally adaptive regression splines [13] (Section 4.1), quantile regression [11] and support vector machines (Section 5).

(b) Our efficient algorithms allow us to pose statistically motivated regularized problems that have not been considered in the literature. In this context, we propose robust versions of the LASSO for regression and classification (Section 3).

**2. The piecewise linear toolbox.** For the coefficient paths to be piecewise linear, we require that $\frac{\partial \hat{\beta}(\lambda)}{\partial \lambda} / \|\frac{\partial \hat{\beta}(\lambda)}{\partial \lambda}\|$ be a piecewise constant vector as a function of $\lambda$. Using Taylor expansions of the normal equations for the minimizing problem (1), we can show that if $L, J$ are both twice differentiable in the neighborhood of a solution $\hat{\beta}(\lambda)$, then

$$(5) \qquad \frac{\partial \hat{\beta}(\lambda)}{\partial \lambda} = -[\nabla^2 L(\hat{\beta}(\lambda)) + \lambda \nabla^2 J(\hat{\beta}(\lambda))]^{-1} \nabla J(\hat{\beta}(\lambda)),$$

where we are using the notation $L(\hat{\beta}(\lambda))$ in the obvious way, that is, we make the dependence on the data $X, \mathsf{y}$ (here assumed constant) implicit.

PROPOSITION 1.  *A sufficient and necessary condition for the solution path to be linear at $\lambda_0$ when $L, J$ are twice differentiable in a neighborhood of $\hat{\beta}(\lambda_0)$ is that*

$$(6) \qquad -[\nabla^2 L(\hat{\beta}(\lambda)) + \lambda \nabla^2 J(\hat{\beta}(\lambda))]^{-1} \nabla J(\hat{\beta}(\lambda))$$

*is a proportional (i.e., constant up to multiplication by a scalar) vector in $\mathbb{R}^p$ as a function of $\lambda$ in a neighborhood of $\lambda_0$.*

Proposition 1 implies sufficient conditions for piecewise linearity:

- $L$ is piecewise quadratic as a function of $\beta$ along the optimal path $\hat{\beta}(\lambda)$, when $X, \mathsf{y}$ are assumed constant at their sample values; and
- $J$ is piecewise linear as a function of $\beta$ along this path.

We devote the rest of this paper to examining some families of regularized problems which comply with these conditions.

On the loss side, this leads us to consider functions $L$ which are:

- Pure quadratic loss functions, like those of linear regression.
- A mixture of quadratic and linear pieces, like Huber's loss [10]. These loss functions are of interest because they generate robust modeling tools. They will be the focus of Section 3.
- Loss functions which are piecewise linear. These include several widely used loss functions, like the hinge loss of the support vector machine (4) and the check loss of quantile regression [11]

$$L(\mathsf{y}, X\beta) = \sum_i l(y_i, \beta^\top \mathbf{x}_i),$$

where

$$(7) \qquad l(y_i, \beta^\top \mathbf{x}_i) = \begin{cases} \tau \cdot (y_i - \beta^\top \mathbf{x}_i), & \text{if } y_i - \beta^\top \mathbf{x}_i \geq 0, \\ (1 - \tau) \cdot (\beta^\top \mathbf{x}_i - y_i), & \text{otherwise,} \end{cases}$$

and $\tau \in (0, 1)$ indicates the quantile of interest.

On the penalty side, our results lead us to consider the $\ell_1$ and $\ell_\infty$ penalties as building blocks for piecewise linear solution paths. For lack of space, we limit the discussion in this paper to the $\ell_1$ penalty and its variants (like total variation penalties discussed in Section 4). Results on the $\ell_\infty$ penalty can be found in a full technical report on the second author's homepage.

$\ell_1$ regularization has several favorable statistical properties. Using $\ell_1$ regularization results in "sparse" solutions with a relatively small fraction of nonzero coefficients, as opposed to $\ell_2$ regularization which forces all nonzero coefficients [17]. In particular, if the number of predictors is larger than the number of observations ($p > n$), then for any $\lambda$ there exists an $\ell_1$-regularized solution with at most $n$ nonzero coefficients [15]. Thus, in situations where the number of relevant variables is small and there are a lot of irrelevant "noise" variables, $\ell_1$ regularization may prove far superior to $\ell_2$ regularization from a prediction error perspective. From an inference/interpretation perspective, $\ell_1$ regularization gives "smooth" variable selection and more compact models than $\ell_2$ regularization. In the case of orthogonal wavelet bases, the soft thresholding method proposed by [2], which is equivalent to $\ell_1$ regularization, is asymptotically nearly optimal (in a minimax sense) over a wide variety of loss functions and estimated functions.

It is not surprising, therefore, that $\ell_1$ regularization and its variants have been widely and successfully used in different fields, including engineering and signal processing (such as basis pursuit and wavelet thresholding), machine learning (such as boosting and $\ell_1$ SVM) and, obviously, statistics, where $\ell_1$ and total variation penalties are prevalent.

**3. Almost quadratic loss functions with $\ell_1$ penalty.** In this section, we first define a family of "almost quadratic" loss functions whose $\ell_1$-penalized versions generate piecewise linear solution paths. We formulate and prove an algorithm, which is an extension of the LAR–LASSO algorithm, that generates the $\ell_1$-penalized solution paths for all members of this family. We then concentrate on two members of this family—Huberized LASSO for regression and $\ell_1$-penalized Huberized squared hinge loss for classification—which define new, robust, efficient and adaptable modeling tools. An R implementation of these tools is available from the second author's homepage, www.stat.lsa.umich.edu/˜jizhu/code/piecewise/.

3.1. *Main results.* We fix the penalty to be the $\ell_1$ penalty,

$$(8) \qquad J(\beta) = \|\beta\|_1 = \sum_j |\beta_j|,$$

and the loss is required to be differentiable and piecewise quadratic in a fixed function of the sample response and the "prediction" $\beta^\top x$,

$$(9) \quad L(\mathbf{y}, X\beta) = \sum_i l(y_i, \beta^\top \mathbf{x}_i), \qquad l(y, \beta^\top \mathbf{x}) = a(r)r^2 + b(r)r + c(r),$$

where $r = (y - \beta^\top \mathbf{x})$ is the residual for regression and $r = (y\beta^\top \mathbf{x})$ is the margin for classification; and $l(r)$ is a quadratic spline, that is, $a(\cdot), b(\cdot), c(\cdot)$ are piecewise constant functions, defined so as to make the function $l$ differentiable.

Some examples from this family are:

- The squared error $l(y, \beta^\top \mathbf{x}) = (y - \beta^\top \mathbf{x})^2$, that is, $a \equiv 1, b \equiv 0, c \equiv 0$.
- Huber's loss function with *fixed* knot $t$,

(10) $$l(y, \beta^\top \mathbf{x}) = \begin{cases} (y - \beta^\top \mathbf{x})^2, & \text{if } |y - \beta^\top \mathbf{x}| \leq t, \\ 2t|y - \beta^\top \mathbf{x}| - t^2, & \text{otherwise.} \end{cases}$$

- Squared hinge loss for classification,

(11) $$l(y, \beta^\top \mathbf{x}) = (1 - y\beta^\top \mathbf{x})_+^2.$$

Note that the hinge loss of the support vector machine and the check loss of quantile regression do not belong to this family as they are not differentiable at $y\beta^\top \mathbf{x} = 1$ and $y - \beta^\top \mathbf{x} = 0$, respectively.

THEOREM 2. *All regularized problems of the form* (1) *using* (8), (9) *(with $r$ being either the residual or the margin) generate piecewise linear optimal coefficient paths $\hat{\beta}(\lambda)$ as the regularization parameter $\lambda$ varies.*

PROOF. We prove the theorem formally using the Karush–Kuhn–Tucker formulation of the optimization problem.

We rewrite the regularized optimization problem as

$$\min_{\beta^+, \beta^-} \sum_i l(y_i, (\beta^+ - \beta^-)^\top \mathbf{x}_i) + \lambda \sum_j (\beta_j^+ + \beta_j^-)$$

$$\text{subject to} \quad \beta_j^+ \geq 0, \beta_j^- \geq 0 \quad \forall j.$$

The Lagrange primal function is

$$\sum_i l(y_i, (\beta^+ - \beta^-)^\top \mathbf{x}_i) + \lambda \sum_j (\beta_j^+ + \beta_j^-) - \sum_j \lambda_j^+ \beta_j^+ - \sum_j \lambda_j^- \beta_j^-.$$

The derivatives of the primal and the corresponding KKT conditions imply

$$(\nabla L(\beta))_j + \lambda - \lambda_j^+ = 0, \qquad -(\nabla L(\beta))_j + \lambda - \lambda_j^- = 0,$$

$$\lambda_j^+ \beta_j^+ = 0, \qquad \lambda_j^- \beta_j^- = 0.$$

Using these we can figure that at the optimal solution for fixed $\lambda$ the following scenarios should hold:

$$\lambda = 0 \Rightarrow (\nabla L(\beta))_j = 0 \qquad \forall j \text{ (unconstrained solution)},$$

$$\beta_j^+ > 0, \lambda > 0 \Rightarrow \lambda_j^+ = 0 \Rightarrow (\nabla L(\beta))_j = -\lambda < 0 \Rightarrow \lambda_j^- > 0 \Rightarrow \beta_j^- = 0,$$

$$\beta_j^- > 0, \lambda > 0 \Rightarrow \beta_j^+ = 0 \qquad \text{(by similar reasoning)},$$

$$|(\nabla L(\beta))_j| > \lambda \Rightarrow \text{contradiction.}$$

Based on these possible scenarios we can see that:

- Variables can have nonzero coefficients only if their "generalized absolute correlation" $|\nabla L(\hat{\beta}(\lambda))_j|$ is equal to $\lambda$. Thus, for every value of $\lambda$ we have a set of "active" variables $\mathcal{A} = \{j : \hat{\beta}_j(\lambda) \neq 0\}$ such that

$$(12) \qquad j \in \mathcal{A} \Rightarrow |\nabla L(\hat{\beta}(\lambda))_j| = \lambda, \qquad \mathrm{sgn}(\nabla L(\hat{\beta}(\lambda))_j) = -\mathrm{sgn}(\hat{\beta}(\lambda)_j),$$

$$(13) \qquad j \notin \mathcal{A} \Rightarrow |\nabla L(\hat{\beta}(\lambda))_j| \leq \lambda.$$

- When $\lambda$ changes, the direction in which $\hat{\beta}(\lambda)$ is moving, that is, $\frac{\partial \hat{\beta}(\lambda)}{\partial \lambda}$, should be such that it maintains the conditions (12), (13).

So, if we know what the "active" set $\mathcal{A}$ is, it is a simple task to check that as long as we are in a region where the loss is twice differentiable and the penalty is right differentiable, we will have

$$(14) \qquad \frac{\partial \hat{\beta}(\lambda)_{\mathcal{A}}}{\partial \lambda} = -(\nabla^2 L(\hat{\beta}(\lambda))_{\mathcal{A}})^{-1} \mathrm{sgn}(\hat{\beta}(\lambda)_{\mathcal{A}}),$$

which is just a version of (5), limited to only the active variables and substituting the $\ell_1$ penalty for $J$.

For the family of almost quadratic loss functions, we can derive $\nabla^2 L(\hat{\beta}(\lambda))_{\mathcal{A}}$ explicitly,

$$\nabla^2 L(\hat{\beta}(\lambda))_{\mathcal{A}} = \sum_i 2a(r(y_i, \hat{\beta}(\lambda)_{\mathcal{A}}^\top \mathbf{x}_{\mathcal{A}i})) \mathbf{x}_{\mathcal{A}i} \mathbf{x}_{\mathcal{A}i}^\top.$$

Since $a(\cdot)$ is a piecewise constant function, then $\nabla^2 L(\hat{\beta}(\lambda))_{\mathcal{A}}$ and $\partial \hat{\beta}(\lambda)_{\mathcal{A}}/\partial \lambda$ are also piecewise constant; therefore, the solution path $\hat{\beta}(\lambda)$ is piecewise linear.

When one of the following "events" occurs, twice differentiability is violated and hence the direction in (14) will change:

- Add a variable: A new variable should join $\mathcal{A}$; that is, we reach a point where $|\nabla L(\hat{\beta}(\lambda))_{\mathcal{A}^C}| \leq \lambda$ will cease to hold if $\hat{\beta}(\lambda)$ keeps moving in the same direction.
- Drop a variable: A coefficient in $\mathcal{A}$ hits 0. In that case, we reach a non-differentiability point in the penalty and we can see that $\mathrm{sgn}(\nabla L(\hat{\beta}(\lambda))_{\mathcal{A}}) = -\mathrm{sgn}(\beta_{\mathcal{A}})$ will cease to hold if we continue in the same direction. Thus we need to drop the coefficient hitting 0 from $\mathcal{A}$.
- Cross a knot: A "generalized residual" $r(y_i, \hat{\beta}(\lambda)^\top \mathbf{x}_i)$ hits a non-twice differentiability point (a "knot") in $L$, for example the "Huberizing" point $t$ in (10), or the hinge point 1 in (11).

So we conclude that the path $\hat{\beta}(\lambda)$ will be piecewise linear, with the direction given by (14) and direction changes occurring whenever one of the three events above happens. When it happens, we need to update $\mathcal{A}$ or $a(r)$ to get a feasible scenario and recalculate the direction using (14). $\quad\square$

Based on the arguments in the proof we can derive a generic algorithm to generate coefficient paths for all members of the "almost quadratic" family of loss functions with $\ell_1$ penalty. The LAR–LASSO algorithm [3] is a simplified version of this algorithm since "knot crossing" events do not occur in the LASSO (as the loss is twice differentiable). Our algorithm starts at $\lambda = \infty$ and follows the linear pieces, while identifying the "events" and recalculating the direction when they occur.

ALGORITHM 1. An algorithm for "almost quadratic" loss with $\ell_1$ penalty.

1. Initialize:

$$\beta = 0, \qquad \mathcal{A} = \arg\max_j |\nabla L(\beta)|_j, \qquad \gamma_{\mathcal{A}} = -\operatorname{sgn}(\nabla L(\beta))_{\mathcal{A}}, \qquad \gamma_{\mathcal{A}^C} = 0.$$

2. While $(\max |\nabla L(\beta)| > 0)$:

   (a) $d_1 = \min\{d > 0 : |\nabla L(\beta + d\gamma)_j| = |\nabla L(\beta + d\gamma)_{\mathcal{A}}|, \ j \notin \mathcal{A}\}$,
       $d_2 = \min\{d > 0 : (\beta + d\gamma)_j = 0, \ j \in \mathcal{A}\}$ (hit 0),
       $d_3 = \min\{d > 0 : r(y_i, (\beta + d\gamma)^\top \mathbf{x}_i)$ hits a "knot," $i = 1, \ldots, n\}$.
       Find step length: $d = \min(d_1, d_2, d_3)$.
   (b) Take step: $\beta \leftarrow \beta + d\gamma$.
   (c) If $d = d_1$ then add variable attaining equality at $d$ to $\mathcal{A}$.
       If $d = d_2$ then remove variable attaining 0 at $d$ from $\mathcal{A}$.
       If $d = d_3$ for $i^*$, then assign new $a(r(y_{i^*}, \beta^\top \mathbf{x}_{i^*}))$ from (9).
   (d) Calculate new direction:

$$C = \sum_i a(r(y_i, \beta^\top \mathbf{x}_i)) \mathbf{x}_{\mathcal{A},i} \mathbf{x}_{\mathcal{A},i}^\top,$$

$$\gamma_{\mathcal{A}} = C^{-1} \cdot \operatorname{sgn}(\beta_{\mathcal{A})} \quad \text{and} \quad \gamma_{\mathcal{A}^C} = 0.$$

It should be noted that our formulation here of the "almost quadratic" family with $\ell_1$ penalty has ignored the existence of a nonpenalized intercept. This has been done for simplicity of exposition, however incorporating a nonpenalized intercept into the algorithm is straightforward.

3.2. *Computational considerations.* What is the computational complexity of running Algorithm 1 on a dataset with $n$ observations and $p$ variables? The major computational cost for each step involves figuring out the step length in (2a), and updating the new direction in (2d). The former takes $O(np)$ calculations, and the latter requires $O(|\mathcal{A}|^2)$ computations by using inverse updating and downdating.

It is difficult to predict the number of steps on the solution path for arbitrary data. According to our experience, the total number of steps taken by the algorithm is on average $O(n)$. This can be heuristically understood as follows. If $n > p$, it takes $O(p)$ steps to add all variables and $O(n)$ steps for knot crossing; if $n < p$, since at most $n$ variables are allowed in the fitted model, it takes $O(n)$ steps for both adding variables and crossing knots; the "drop events" are usually rare, $O(1)$. Since the maximum value of $|\mathcal{A}|$ is $\min(n, p)$, it suggests the overall computational cost is $O(n^2 p)$.

3.3. *The Huberized LASSO for regression.* We now concentrate on two members of the "almost quadratic" family of loss functions—one for regression and one for classification.

We first consider the Huberized LASSO for regression. The loss is given by (10). It is robust in the sense defined in [10], in that it protects against "contamination" of the assumed normal errors. It is "almost quadratic" as defined in Section 3.1, and so Theorem 2 and Algorithm 1 apply to its $\ell_1$ regularized solution paths.

*Prostate cancer dataset.* We use the "prostate cancer" dataset [17] to compare the prediction performance of the Huberized LASSO to that of the LASSO on the original data and after we artificially "contaminate" the data by adding large constants to a small number of responses.

We used the training-test split as in [8]. The training set consists of 67 observations and the test set of 30 observations. We ran the LASSO and the Huberized LASSO with a knot at $t = 1$ on the original dataset, and on the "contaminated" dataset where 5 has been added/subtracted to the responses of 12 observations.

Figure 1 shows the mean squared error on the 30 test set observations for the four resulting regularized solution paths from solving the LASSO and Huberized LASSO for all possible values of $\lambda$ on the two datasets. We observe that on the noncontaminated data, the LASSO (solid) and Huberized LASSO (dashed) perform quite similarly. When we add contamination, the Huberized LASSO (dash-dotted) does not seem to suffer from it at all, in that its best test set performance is comparable to that of both regularized models on the noncontaminated data. The prediction performance of the standard LASSO (dotted), on the other hand, deteriorates significantly ($t$-test $p$-value 0.045) when contamination is added, illustrating the lack of robustness of squared error loss.

The two LASSO solutions contain nine linear pieces each, while the Huber-LASSO path for the noncontaminated data contains 41 pieces, and the one for the contaminated data contains 39 pieces; both agree with our conjecture in Section 3.2 that the number of steps is $O(n)$. Figure 2 shows the solution paths for the contaminated LASSO model and the contaminated Huber-LASSO model. We observe that the two paths are quite different and the two best models (corresponding to the solid vertical lines) are also different.

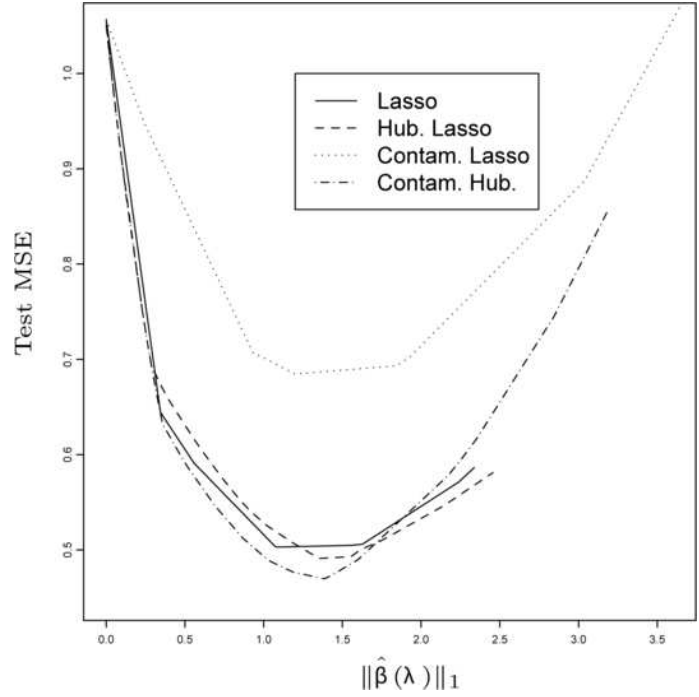

FIG. 1.  *Test MSE of the models along the regularized paths. See text for details.*

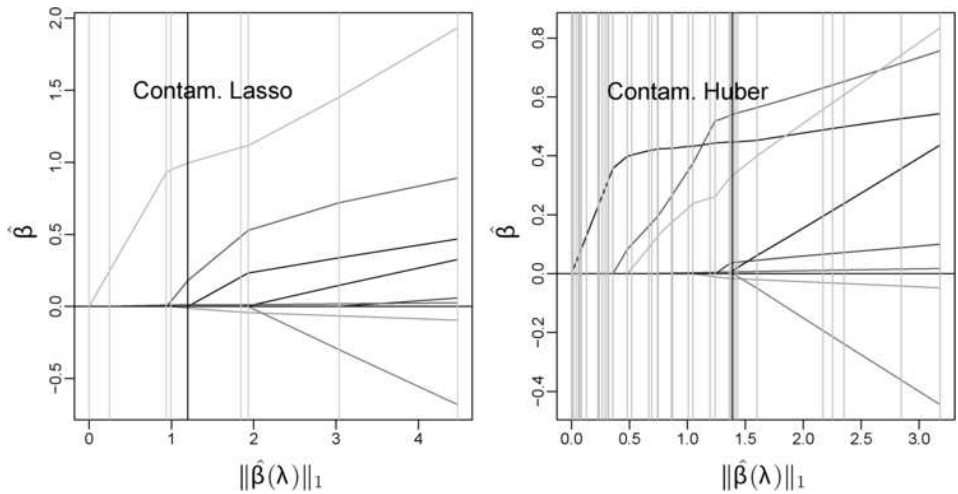

FIG. 2.  *Solution paths of the LASSO (left) and the Huberized LASSO (right) on the contaminated prostate cancer training data. The vertical grey lines correspond to the steps along the solution paths. The vertical solid lines correspond to the models that give the best performances on the test data.*

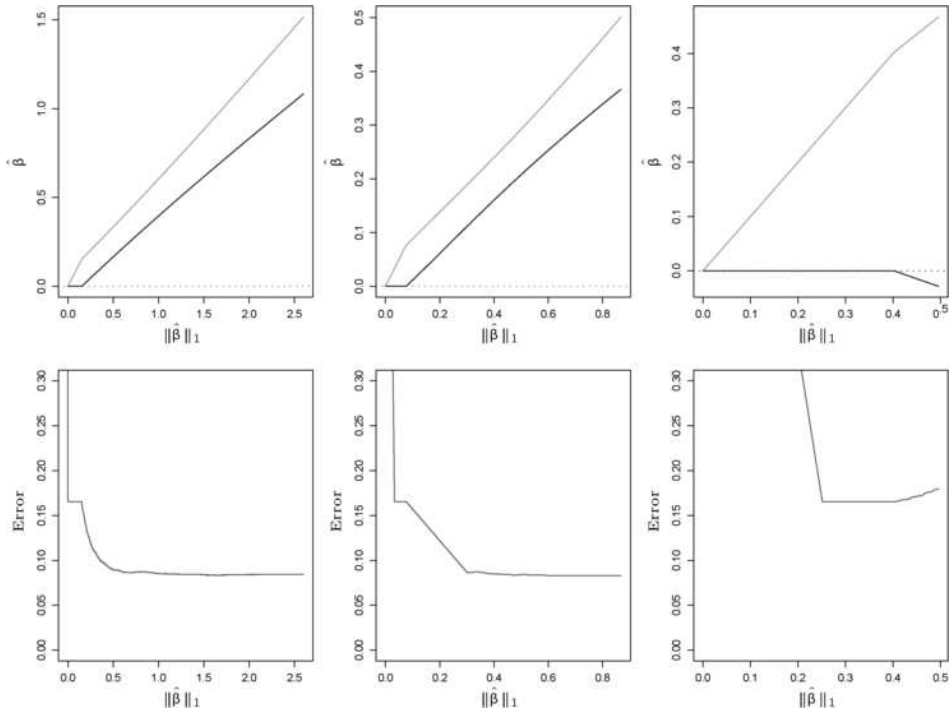

Fig. 3. *Regularized paths and prediction errors for the logistic loss (left), Huberized squared hinge loss (middle) and squared hinge loss (right). The logistic loss and the Huberized squared hinge loss are both less affected by the outlier.*

3.4. *The Huberized squared hinge loss for classification.* For classification we would like to have a loss which is a function of the margin, $r(y, \beta^\top \mathbf{x}) = (y\beta^\top \mathbf{x})$. This is true of all loss functions typically used for classification, like the negative binomial log-likelihood for logistic regression, the hinge loss for the support vector machine and exponential loss for boosting. The properties we would like from our classification loss are:

- We would like it to be "almost quadratic," so we can apply the Algorithm 1 in Section 3.1.
- We would like it to be robust, that is, linear for large absolute value negative margins (like the logistic or hinge), so that outliers will have a small effect on the fit.

This leads us to suggest for classification the "Huberized squared hinge loss," that is, (11) "Huberized" at $t < 1$,

$$(15) \quad l(y, \beta^\top \mathbf{x}) = \begin{cases} (1-t)^2 + 2(1-t)(t - y\beta^\top \mathbf{x}), & \text{if } y\beta^\top \mathbf{x} \le t, \\ (1 - y\beta^\top \mathbf{x})^2, & \text{if } t < y\beta^\top \mathbf{x} \le 1, \\ 0, & \text{otherwise.} \end{cases}$$

It is a simple task to show that

$$\arg\min_f \mathbf{E}_y l(y, f) = 2\Pr(y = 1) - 1.$$

Hence the population minimizer of the Huberized squared hinge loss gives the correct sign for classification.

To illustrate the robustness of this loss (15) and its computational superiority over the logistic loss, we considered the following simple example: $\mathbf{x} \in \mathbb{R}^2$ with class centers at $(-1, -1)$ (class "−1") and $(1, 1)$ (class "1") with one big outlier at $(30, 100)$ belonging to the class "−1." The Bayes model, ignoring the outlier, is to classify to class "1" if and only if $x_1 + x_2 > 0$.

Figure 3 shows the regularized solution paths and misclassification rate for this example using the logistic loss (left), the Huberized squared hinge loss (middle) and the squared hinge loss (right), all with $\ell_1$ penalty. We observe that the logistic and Huberized regularized model paths are both less affected by the outlier than the non-Huberized squared loss. However, logistic loss does not allow for efficient calculation of the $\ell_1$ regularized path.

## 4. Nonparametric regression, total variation penalties and piecewise linearity.
Total variation penalties and closely associated spline methods for nonparametric regression have experienced a surge of interest in the statistics literature in recent years. The total variation of a univariate differentiable function $f(x)$ is

$$TV_{\text{dif}}(f) = \int_{-\infty}^{\infty} |f'(x)| \, dx.$$

If $f$ is nondifferentiable on a countable set $x_1, x_2, \ldots$, then $TV(f)$ is the sum of $TV_{\text{dif}}(f)$, calculated over the differentiable set only and the absolute "jumps" in $f$ where it is noncontinuous. In what follows we assume the range of $f$ is limited to $[0, 1]$.

Total variation penalties tend to lead to regularized solutions which are polynomial splines. [13] investigates the solutions to total-variation penalized least squares problems. The authors use total variation of $(k-1)$st order derivatives,

$$(16) \qquad \sum_{i=1}^{n} (y_i - f(x_i))^2 + \lambda \cdot TV(f^{(k-1)}).$$

They show that there always exists a solution $\hat{f}_{k,\lambda}$ such that $\hat{f}_{k,\lambda}^{(k-1)}$ is piecewise constant, that is, $\hat{f}_{k,\lambda}$ is a polynomial spline of order $k$. For $k \in \{1, 2\}$ the knots of the spline solutions are guaranteed to be at the data points.

A similar setup is considered in [1]. Their *taut-string* and *local squeezing* methods lead to solutions that are polynomial splines of degree 0 or 1, with knots at data points.

Now, consider a polynomial spline $f$ of order $k$, with $h$ knots located at $0 < t_1 < \cdots < t_h < 1$, that is,

$$(17) \qquad f(x) = \sum_{j=1}^{h} \beta_j (x - t_j)_+^{k-1} + q(x),$$

where $q(x)$ is a polynomial of degree $k-1$. The total variation of the $(k-1)$st derivative of $f$ clearly corresponds to an $\ell_1$ norm of the set of coefficients of the appropriate spline basis functions, anchored at the knots,

$$(18) \qquad TV(f^{(k-1)}) = (k-1)! \cdot \sum_{j=1}^{h} |\beta_j|.$$

If the knots $t_1, \ldots, t_h$ are fixed in advance (e.g., at the data points), then a total variation penalized problem is equivalent to an $\ell_1$-penalized regression problem, with $p = h$ derived predictors. If we also employ squared error loss, we get a LASSO problem, and we can use the LAR–LASSO algorithm to compute the complete regularized solution path. The only difference from the standard LASSO is the existence of $k$ nonpenalized coefficients for the polynomial $q(x)$, instead of the intercept only for the LASSO. This requires only a slight modification to the LAR–LASSO algorithm. This leads to essentially the same algorithm as Algorithm 2 of [13] for finding the regularized path for any $k$ with a *fixed, predetermined* set of candidate knots.

4.1. *Locally adaptive regression splines.* We now concentrate on the family of penalized problems (16) defined by Mammen and van de Geer [13]. As we have mentioned, [13] develops an exact method for finding $\hat{f}_{k,\lambda}$ when $k \in \{1, 2\}$ and *approximate* methods for $k > 2$ (where the knots of the optimal solutions are not guaranteed to be at the data points). We now show how we can use our approach to find the spline solution $\hat{f}_{k,\lambda}$ *exactly* for any natural $k$. The resulting algorithms get practically more complicated as $k$ increases, but their theoretical computational complexity remains fixed.

When the knots are not guaranteed to be at the data points, we can still write the total variation of polynomial splines as the sum of $\ell_1$ norms of coefficients of basis functions, as in (18). However, we do not have a finite predefined set of candidate basis functions. Rather, we are dealing with an infinite set of candidate basis functions of the form

$$\mathcal{X} = \{(x - t)_+^{k-1} : 0 \le t \le 1\}.$$

Our algorithm for tracking the regularized solution path $\hat{f}_{k,\lambda}$ in this case proceeds as follows. We start at the solution for $\lambda = \infty$, which is the least squares $(k-1)$st degree polynomial fit to the data. Given a solution $\hat{f}_{k,\lambda_0}$ for some value of $\lambda_0$, which includes $n_{\lambda_0}$ knots at $t_1, \ldots, t_{n_{\lambda_0}}$, denote

$$\mathbf{z}(x) = (1, x, x^2, \ldots, x^{k-1}, (x - t_1)_+^{k-1}, \ldots, (x - t_{n_{\lambda_0}})_+^{k-1})^\top,$$

which is the current predictor vector. Following (17) we can write

$$\hat{f}_{k,\lambda_0}(x) = \hat{\beta}(\lambda_0)^\top \mathbf{z}(x).$$

Following the logic of the LAR–LASSO algorithm, we see that the solution will change as

$$\hat{f}_{k,\lambda_0-d} = (\hat{\beta}(\lambda_0) + d\gamma)^\top \mathbf{z}(x),$$

where $\gamma = -(k-1)! \cdot (Z^\top Z)^{-1} \cdot \mathbf{s}$, $Z = (\mathbf{z}(x_1), \dots, \mathbf{z}(x_n))^\top$ and $\mathbf{s} \in \mathbb{R}^{k+n_{\lambda_0}}$ is a vector with 0 components corresponding to $1, x, \dots, x^{k-1}$ and $\pm 1$ components corresponding to each $(x-t_j)_+^{k-1}$ [with the sign being the opposite of the sign of $(x-t_j)_+^{k-1\top}(\mathbf{y} - \hat{f}_{k,\lambda_0})$]. What we now need to identify is the value of $\lambda$ at which an additional knot needs to be added, and the location of that knot. Consider first a fixed knot candidate $t$. Then we can see that the LAR–LASSO criterion for adding this knot to the set of "active" knots is

$$|\mathbf{x}_t^\top(\mathbf{y} - Z\hat{\beta}(\lambda_0) - (\lambda_0 - \lambda)Z\gamma)| = \lambda,$$

where $\mathbf{x}_t = (\mathbf{x} - t)_+^{k-1}$ (column vector of length $n$). More explicitly, define

$$(19) \qquad \lambda_+(t) = \frac{\mathbf{x}_t^\top(\mathbf{y} - Z\hat{\beta}(\lambda_0) - \lambda_0 Z\gamma)}{1 - \mathbf{x}_t^\top Z\gamma},$$

$$(20) \qquad \lambda_-(t) = \frac{\mathbf{x}_t^\top(\mathbf{y} - Z\hat{\beta}(\lambda_0) - \lambda_0 Z\gamma)}{-1 - \mathbf{x}_t^\top Z\gamma}.$$

Then we can write

$$(21) \qquad \lambda(t) = \begin{cases} \max(\lambda_+(t), \lambda_-(t)), & \text{if } \max(\lambda_+(t), \lambda_-(t)) \le \lambda_0, \\ \min(\lambda_+(t), \lambda_-(t)), & \text{if } \max(\lambda_+(t), \lambda_-(t)) > \lambda_0. \end{cases}$$

Now we see that we can in fact let $t$ be a parameter and find the next knot to be added to the optimal solution path by maximizing $\lambda(t)$, that is,

$$(22) \qquad \lambda_{\text{add}} = \max_{t \in (0,1) \setminus \{t_1, \dots, t_{n_{\lambda_0}}\}} \lambda(t),$$

which is the value of $\lambda$ where we stop moving in direction $\gamma$, add a knot at the argument of the maximum and recalculate the direction $\gamma$.

Solving (22) requires finding the local extrema of the functions in (21), which are rational functions within each interval between two points which are either data points or knots (with numerator and denominator both of degree $k-1$). Thus, a reasonable tactic is to find the extrema within each such interval, then compare them between the intervals to find the overall solution to (22). For smaller values of $k$ it can be solved manually and exactly:

- For $k \in \{1, 2\}$, we get a ratio of constant or linear functions in (21), and therefore the extrema—and the knots—are guaranteed to be at the data points, leading to the algorithm of [13].
- For $k = 3$ we get a ratio of quadratics in (21), and we can find the extrema within each segment analytically. These extrema may not correspond to the segment's end points, and so we may have knots that are not at data points.

Assuming we have the code to solve the maximization problem in (22), Algorithm 2 gives a general schema for following the solution path $\hat{f}_{k,\lambda}$ for any value of $k$.

ALGORITHM 2.  Tracking the path of TV-penalized solutions.

1. Initialize:

   $f(x) = (1, x, \ldots, x^{k-1})^\top \beta_{\mathrm{ls}}$ is the LS polynomial fit of degree $k - 1$,

   $u = \underset{t \in (0,1)}{\arg\max} |(\mathbf{x} - t)_+^{k-1\top} (\mathbf{y} - f(\mathbf{x}))|$        (assumed unique),

   $\mathcal{T} = \{u\}, \qquad \lambda_0 = (k-1)! \cdot |(\mathbf{x} - u)_+^{k-1\top} (\mathbf{y} - f(\mathbf{x}))|,$

   $Z = (\mathbf{1}, \mathbf{x}, \ldots, \mathbf{x}^{k-1}, (\mathbf{x} - u)_+^{k-1}), \qquad \hat{\beta}(\lambda_0) = (\beta_{\mathrm{ls}}^\top, 0)^\top,$

   $\mathbf{s} = (\mathbf{0}_k^\top, -\operatorname{sgn}\{(\mathbf{x} - u)_+^{k-1\top} (\mathbf{y} - f(\mathbf{x}))\})^\top.$

2. While $\sum_i (y_i - f(x_i))^2 > 0$:

   (a) Set $\gamma = -(k-1)!(Z^\top Z)^{-1} \mathbf{s}$.

   (b) $\forall t \in (0,1) \setminus \mathcal{T}$ define $\lambda_+(t), \lambda_-(t), \lambda(t)$ as in (19)–(21).

   (c) Solve the maximum problem in (22) to get $\lambda_{\mathrm{add}}$.

   (d) Let $\lambda_{\mathrm{rem}} = \lambda_0 - \min\{d > 0 : \exists j > k \text{ s.t. } \hat{\beta}_j(\lambda_0) + d\gamma_j = 0\}$.

   (e) If $\lambda_{\mathrm{add}} > \lambda_{\mathrm{rem}}$ add a knot at the point attaining the maximum in (22), and update $\mathcal{T}$, $Z$ and $\mathbf{s}$.

   (f) Similarly, if $\lambda_{\mathrm{add}} < \lambda_{\mathrm{rem}}$ remove the knot attaining 0 at $\lambda_{\mathrm{rem}}$.

   (g) In both cases, update:

   $$\hat{\beta}(\lambda_0) \leftarrow \hat{\beta}(\lambda_0) + (\lambda_0 - \max(\lambda_{\mathrm{add}}, \lambda_{\mathrm{rem}}))\gamma,$$

   $$\lambda_0 \leftarrow \max(\lambda_{\mathrm{add}}, \lambda_{\mathrm{rem}}).$$

Since we can never have more than $n - k$ knots in a solution $\hat{f}_{k,\lambda}$ [15], the computational complexity of each iteration of the algorithm is bounded at $O(n^2)$ calculations for finding the next knot and $O(n^2)$ for calculating the next direction (using updating formulas). The number of steps of the algorithm is difficult to bound, but from our experience seems to behave like $O(n)$ (which is the conjecture of [13]).

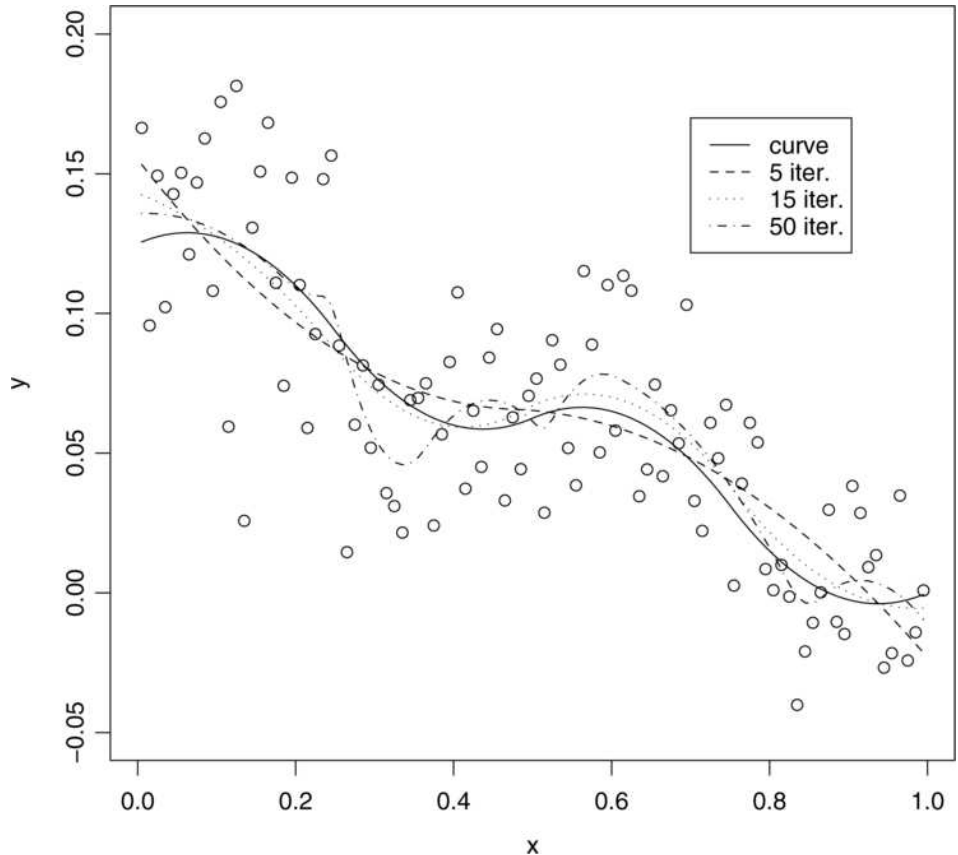

Fig. 4.  *Applying Algorithm 2 (with $k = 3$) to a data example where the underlying curve is a quadratic spline with knots at 0.25, 0.5 and 0.75. See text for details.*

4.2. *Simple data example*: $k = 3$.  We illustrate our algorithm on a simple data example. We select 100 $x$ samples uniformly on $(0, 1)$. We draw the corresponding $y$ values as $N(g(x), 0.03^2)$, where $g(x)$ is a polynomial spline with knots at $0.25, 0.5$ and $0.75$,

$$g(x) = 0.125 + 0.125x - x^2 + 2(x - 0.25)_+^2 - 2(x - 0.5)_+^2 + 2(x - 0.75)_+^2.$$

$g(x)$ is plotted as the solid line in Figure 4, and the noisy $y$ values as circles. The signal-to-noise ratio is about 1.4.

We apply our Algorithm 2 with $k = 3$. Figure 4 shows the resulting models after 5, 15 and 50 iterations of the algorithm. After 5 iterations, the regularized spline contains three knots like the true $g$, but these are all around 0.5. The fitted model, drawn as a dashed curve, is clearly underfitted. The corresponding reducible squared prediction error is $9.5 \times 10^{-4}$.

After 15 iterations, the spline contains four knots, at 0.225, 0.255, 0.485 and 0.755. The first one has a small coefficient, and the other three closely correspond to the knots in $g$. The resulting fit (dotted curve) is a reasonable approximation of $g$, and the reducible squared error is about $3.1 \times 10^{-4}$.

After 50 iterations the model contains ten knots and the data is clearly overfitted (dash-dotted curve, reducible squared error $8.2 \times 10^{-4}$).

Although the algorithm should in principle continue until it interpolates the data, in practice it terminates before (in this case after about 180 iterations) and is numerically unable to further improve the fit. This is analogous to the situation described in [7] for kernel SVM, where the effective rank of the kernel matrix is significantly smaller than $n$, since many eigenvalues are effectively zero.

**5. Using $\ell_1$ loss and its variants.** Piecewise linear nondifferentiable loss functions appear in practice in both regression and classification problems. For regression, absolute value loss variants like the quantile regression loss are quite popular [11]. For classification, the hinge loss is of great importance, as it is the loss underlying the support vector machine [20]. Here we consider a generalized formulation, which covers both of (7) and (4). The loss function has the form

$$(23) \qquad l(r) = \begin{cases} b_1 \cdot |a + r|, & \text{if } a + r \geq 0, \\ b_2 \cdot |a + r|, & \text{if } a + r < 0, \end{cases}$$

with the generalized "residual" being $r = (y - \beta^\top \mathbf{x})$ for regression and $r = (y \cdot \beta^\top \mathbf{x})$ for classification.

When these loss functions are combined with $\ell_1$ penalty (or total variation penalty, in appropriate function classes [12]), the resulting regularized problems can be formulated as linear programming problems. When the path of regularized solutions $\hat{\beta}(\lambda)$ is considered, it turns out to have interesting structure with regard to $\lambda$:

PROPOSITION 3. *For loss functions of the form* (23), *there exists a set of values of the regularization parameter* $0 < \lambda_1 < \cdots < \lambda_m = \infty$ *such that*:

- *The solution $\hat{\beta}(\lambda_k)$ is not uniquely defined, and the set of optimal solutions for each $\lambda_k$ is a straight line in $\mathbb{R}^p$.*
- *For any $\lambda \in (\lambda_k, \lambda_{k+1})$, the solution $\hat{\beta}(\lambda)$ is fixed and equal to the minimum $\ell_1$ norm solution for $\lambda_k$ and the maximum $\ell_1$ norm solution for $\lambda_{k+1}$.*

Proposition 3 generalizes observations on the path of solutions made in the context of quantile regression in [12] and in the context of 1-norm support vector machines in [21]. Note that this leads to describing a regularized path which is *piecewise constant* as a function of the regularization parameter $\lambda$, with jumps at the values $\lambda_1, \ldots, \lambda_m$. However, it is still *piecewise linear* in the $\ell_1$ norm of the solution, $\|\hat{\beta}(\lambda)\|_1$. The algorithm for computing the solution path follows the spirit of our earlier work [21]. For brevity we omit the details. We note, however, that it is fundamentally *different* from the LARS–LASSO algorithm and Algorithm 1, because we are now dealing with a nondifferentiable loss function.

An interesting variant of piecewise linear loss is to replace the $\ell_1$ loss with an $\ell_\infty$ loss, which is also piecewise linear and nondifferentiable. It leads to interesting "mini-max" estimation procedures, popular in many areas, including engineering and control. For example, [19] proposes the use of $\ell_1$-penalized $\ell_\infty$-loss solutions in an image reconstruction problem (but does not consider the solution path). Path-following algorithms can be designed in the same spirit as the $\ell_1$ loss case.

**6. Conclusion.** In this paper we combine computational and statistical considerations in designing regularized modeling tools. We emphasize the importance of both appropriate regularization and robust loss functions for successful practical modeling of data. From a statistical perspective, we can consider robustness and regularization as almost independent desirable properties dealing with different issues in predictive modeling:

- Robustness mainly protects us against wrong assumptions about our error model. It does little or nothing to protect us against the uncertainty about our model structure which is inherent in the finiteness of our data. For example, if our errors really are normal, then squared error loss minimizes the asymptotic variance of the coefficients, no matter how little data we have or how inappropriate our model is [10]. Using robust loss in such a situation is *always* counter-productive.
- Regularization deals mainly with the uncertainty about our predictive model structure by limiting the model space. Note, in this context, the equivalence between the "penalized" formulation (1) and a "constrained" formulation $\min_\beta L(\mathbf{y}, X\beta)$ subject to $J(\beta) \leq s$. The two formulations share the same solution path. The constrained formulation exposes the goal of regularization as "simplifying" the model estimation problem by limiting the set of considered models.

There are many interesting directions in which our work can be extended. We may ask, how can our geometric understanding of the regularized solution paths help us to analyze the statistical properties of the models along the path? For example, [3] has offered analysis of the LASSO path. This becomes much more challenging once we stray away from squared error loss. We may also consider more complex penalty structure, such as local or data-dependent penalties [1] or multiple penalties [18].

Finally, it is worth noting that limiting our discussion to convex problems, for which efficient algorithms can be designed, leaves out some other statistically well motivated fitting approaches. The use of a nonconvex penalty was advocated by Fan and collaborators in several papers [4, 5]. They expose the favorable variable selection property of the penalty function they offer, which can be viewed as an improvement over the use of $\ell_1$ penalty. [16] advocates the use of nonconvex $\psi$-loss in the classification setting, minimizing the effect of outliers and misclassified points. This approach can be viewed as an even more robust version of our Huberized loss function, with strong statistical motivation in terms of asymptotic behavior.

**Acknowledgments.** We thank the referee, Associate Editor and Co-Editors J. Fan and J. Marden for their thoughtful and useful comments, and in particular for introducing us to the relevant literature on total variation penalties. We thank B. Efron, J. Friedman, T. Hastie, R. Tibshirani, B. Yu and T. Zhang for their helpful comments and suggestions.

PREDICTIVE MODELING GROUP                     DEPARTMENT OF STATISTICS
IBM T. J. WATSON RESEARCH CENTER              UNIVERSITY OF MICHIGAN
YORKTOWN HEIGHTS, NEW YORK 10598              1085 SOUTH UNIVERSITY
USA                                          ANN ARBOR, MICHIGAN 48109
E-MAIL: srosset@us.ibm.com                   USA
                                             E-MAIL: jizhu@umich.edu